\documentclass[a4paper,11pt]{article}

\usepackage{amsmath}
\usepackage{amssymb}
\usepackage{amsthm}
\usepackage{amsfonts}
\usepackage{mathtools}
\usepackage{ascmac}
\usepackage{fancybox} 
\usepackage[mathscr]{eucal}
\usepackage{footnpag}
\usepackage{siunitx}
\usepackage{multicol}
\usepackage{bbm}
\usepackage{here}

\usepackage{graphicx}
\usepackage{xcolor}
\usepackage{wrapfig}
\usepackage{float}

\theoremstyle{definition}
\newtheorem{theorem}{Theorem}[section]
\newtheorem{proposition}{Proposition}[section]
\newtheorem{lemma}{Lemma}[section]
\newtheorem{corollary}{Corollary}[section]

\title{ Fibonacci and Lucas numbers arising from \\ two-component spanning forests of wheel graphs}
\author{Tsuyoshi Miezaki \thanks{Faculty of Science and Engineering, Waseda University, Tokyo, 169-8555, Japan, e-mail: miezaki@waseda.jp }, 
Shunya Tamura\thanks{Okegawa City Okegawa West Junior High School, Saitama, 363-0027, Japan, e-mail: shunya.tamura059@gmail.com}
}

\begin{document}
\maketitle

\begin{abstract}
In this paper, we present a constructive bijection between a conditioned spanning forest of the wheel graph $W_{n+1}$ and a spanning tree of the fan graph $F_n$. In addition, by applying the effective resistance formula obtained by Bapat and Gupta \cite{bapat-gupta}, we derive an explicit formula for the number of two-component spanning forests of $W_{n+1}$ in which two specified vertices $u$ and $v$ lie in distinct components. 
Based on this result, we obtain explicit formulas for the following three conditioned two-component spanning forests $F_{W_{n+1}}(v_1\mid v_2)$, $F_{W_{n+1}}(v_1\mid v_3)$, and $F_{W_{n+1}}(v_1\mid v_c)$. These formulas are 
$F_{W_{n+1}}(v_1\mid v_2)=2(f_{2n-1}-1)$, $F_{W_{n+1}}(v_1\mid v_3)=2(\ell_{2n-2}-3)$, $F_{W_{n+1}}(v_1\mid v_c)=f_{2n}$, 
where $f_i$ and $\ell_j$ denote the $i$-th Fibonacci number and $j$-th Lucas number, respectively. 
As these identities show, the enumerations naturally lead to formulas involving Fibonacci numbers and Lucas numbers. 
Taken together, these two approaches show a unified perspective. One is the constructive combinatorial bijection, and the other is the analytic method based on effective resistance. Together they provide a new integrated framework for studying the structure of spanning forests on $W_{n+1}$. 

\end{abstract}



\section{Introduction}

We consider the wheel graph $W_{n+1}$, which consists of one central vertex and an $n$-vertex cycle. 
The number of spanning trees of $W_{n+1}$ is known to be $\ell_{2n}-2$, where $\ell_{j}$ denotes the Lucas numbers. 
This classical result goes back to Myers \cite{myers} and Bogdanowicz \cite{P1}. 
Since Lucas numbers are closely related to Fibonacci numbers, this formula is a typical and well-known example illustrating the combinatorial structure of $W_{n+1}$. 

On the other hand, Bapat and Gupta \cite{bapat-gupta} derived explicit formulas for the effective resistance of wheel graphs and fan graphs, which in turn are connected to the number of two-component spanning forests. 
Recently, Tamura and Tanaka \cite{tamura-tanaka} computed the number of conditioned two-component spanning forests via hitting times between all pairs of vertices on $W_{n+1}$. 

In this paper, we present a constructive bijection between two-component spanning forests of $W_{n+1}$ in which one component contains the central vertex and the other consists of cycle vertices and the spanning trees of the fan graph $F_n$. 
Moreover, we determine the number of such spanning forests. 

Furthermore, by applying the effective resistance formula obtained by Bapat and Gupta \cite{bapat-gupta}, we derive an explicit formula for the number of two-component spanning forests of $W_{n+1}$ in which two specified vertices $u$ and $v$ lie in distinct components. 
As a corollary, we obtain explicit formulas for three types of conditioned two-component spanning forests. 

The rest of this paper is organized as follows.
In Section \ref{sec02}, we define the graphs under consideration and introduce conditioned two-component spanning forests, together with our main theorem. 
In Section \ref{sec03}, we prove the bijection between the two-component spanning forests of $W_{n+1}$ and the spanning trees of $F_n$. 
In Section \ref{sec04}, we provide the proof of the formula for the number of two-component spanning forests using effective resistance. 
Finally, Section \ref{sec05} summarizes our results. 

\section{Definition and main theorem \label{sec02}}

In this section, we define the graphs under consideration, as well as spanning trees, spanning forests, and the fundamental terminology. 
Using these definitions, we describe our main theorem. 

The wheel graph $W_{n+1}$ has a central vertex $v_c$ and an $n$-vertex cycle. 
The central vertex is connected to each cycle vertex $v_i$ ($i=1, 2, \dots, n$). 
We call the edge $\{v_c, v_i\}$ ($i=1, 2, \dots, n$) a spoke, and the edge $\{v_i, v_{i+1}\} \cup \{v_n, v_1\}$ ($i=1, 2, \dots, n-1$) a cycle edge. 

Next, the fan graph $F_m$ consists of a hub vertex $v_h$ and a path on vertices $u_1, u_2, \dots, u_m$. 
The hub vertex $v_h$ is connected to all vertices $u_i$ ($i=1, 2, \dots, m$). 
We call the edges $\{v_h, u_i\}$ hub edges, and the edges $\{u_i, u_{i+1}\}$ path edges. 
In this paper, we treat the fan graph $F_n$. 

As a basic graph-theoretic notion, we consider spanning trees of a graph $G$. 
A spanning tree of $G$ is a tree containing all vertices of $G$. 
We denote the set of spanning trees by $\mathscr{T}(G)$, and its cardinality by $|\mathscr{T}(G)|$. 
In particular, for the fan graph $F_n$, we often write $\mathscr{T}_n := \mathscr{T}(F_n)$. 
Moreover, a forest of $G$ is a subgraph containing all vertices of $G$ and containing no cycles. 

In this paper, let $F_{W_{n+1}}(v_x \mid v_y)$ denote the number of two-component spanning forests in which the two vertices $v_x$ and $v_y$ belong to different components. 
In $W_{n+1}$, the set of spanning forests consisting of a tree containing the central vertex $v_c$ and a tree consisting only of cycle vertices is denoted by $\tau_{n+1}$. 
For each $T \in \tau_{n+1}$, we denote the component containing the central vertex by $T_{\mathrm{center}}$, and the component consisting only of cycle vertices by $T_{\mathrm{cycle}}$. 
In addition, the number of spanning trees of the wheel graph $W_{n+1}$ is known; Myers \cite{myers} and Bogdanowicz \cite{P1} investigated in detail the relation between Fibonacci numbers and Lucas numbers. 
The Fibonacci numbers $f_i$ satisfy the recurrence relation $f_{i+2}=f_{i+1}+f_i$, $f_0=0$, $f_1=1$. 
As a similar sequence, the Lucas numbers $\ell_j$ satisfy $\ell_{j+2}=\ell_{j+1}+\ell_j$, $\ell_0=2$,  $\ell_1=1$. 
Under the above notation, we describe the main theorem.

\begin{theorem} \label{thm01}
The set of two-component spanning forests $\tau_{n+1}$ of $W_{n+1}$ 
is in bijection with the set of spanning trees $\mathscr{T}_{n}$ of $F_n$. 
That is, $\tau_{n+1} \longleftrightarrow \mathscr{T}_{n}$. 
\end{theorem}

For example, 
\[
T=\{v_1v_2,\ v_2v_3,\ v_c v_4\}\in\tau_{4}, \qquad
S=\{u_1u_2,\ u_2u_3,\ v_hu_3\}\in\mathscr{T}_{3}
\]
correspond to each other. 

\begin{figure}[H]
\centering
\includegraphics[width=8cm]{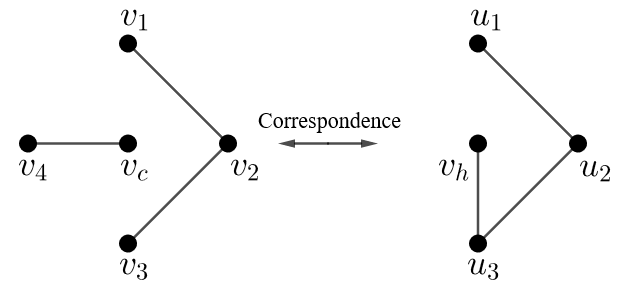}
\caption{Example of the bijection}
\end{figure}

Next, we explain the number of two-component spanning forests of $W_{n+1}$ in which two vertices $u$ and $v$ lie in different components. 

\begin{theorem} \label{thm02}
Let $u$ and $v$ be two vertices of $W_{n+1}$, and let
\[
k=\min(|u-v|,\; n-|u-v|)
\]
be the cycle distance between them. 
Then the number of two-component spanning forests of $W_{n+1}$ in which $u$ and $v$ belong to distinct components is given by
\begin{align*}
F_{W_{n+1}}(u \mid v)
&= f_{2k}\,\mathscr{T}(W_{n+1}) - f_{2n}\,\mathscr{T}(W_{k+1}) \\
&= f_{2k}(\ell_{2n}-2)\;-\;f_{2n}(\ell_{2k}-2).
\end{align*}

Here, $f_i$ denotes the $i$th Fibonacci number ($f_0=0$, $f_1=1$), 
and $\ell_j$ denotes the $j$th Lucas number ($\ell_0=2$, $\ell_1=1$). 
Moreover, $\mathscr{T}(W_k)$ denotes the number of spanning trees of $W_k$. 
\end{theorem}

We apply Theorem \ref{thm02} to $W_{n+1}$ and obtain the number of three types of
two-component spanning forests, whose formulas involve Fibonacci and Lucas numbers. 

\begin{proposition} \label{prop1a}
The numbers of the three types of two-component spanning forests, 
$F_{W_{n+1}}(v_1\mid v_2)$, $F_{W_{n+1}}(v_1\mid v_3)$, and $F_{W_{n+1}}(v_1\mid v_c)$, 
are given by
\begin{itemize}
\item[$\bullet$]
$F_{W_{n+1}}(v_1\mid v_2)=2(f_{2n-1}-1)$

\item[$\bullet$] 
$F_{W_{n+1}}(v_1\mid v_3)=2(\ell_{2n-2}-3)=2(\mathscr{T}(W_{n})-1)$

\item[$\bullet$]
$F_{W_{n+1}}(v_1\mid v_c)=f_{2n}$

\end{itemize} 
Here, $f_i$ denotes the $i$th Fibonacci number ($f_0=0$, $f_1=1$),
and $\ell_j$ denotes the $j$th Lucas number ($\ell_0=2$, $\ell_1=1$).
Moreover, $\mathscr{T}(W_i)$ denotes the number of spanning trees of $W_i$. 
\end{proposition}

\section{Proof of Theorem \ref{thm01}} \label{sec03}
We complete the proof of Theorem \ref{thm01} by using Lemmas \ref{lem01} to \ref{lem05}. 

\begin{lemma} \label{lem01}
For every $T \in \tau_{n+1}$, the cycle component $T_{\mathrm{cycle}}$ is a path consisting of a consecutive segment of the cycle. 
\end{lemma}

\noindent
{\bf Proof. }
Since $T_{\mathrm{cycle}}$ is a connected subgraph of $C_n$ and contains no cycles, it must be a path. 
\qed

\begin{lemma} \label{lem02}
For every $T_{\mathrm{cycle}} \in \tau_{n+1}$, by rotating the labels of the cycle, we may express it as
\[
T_{\mathrm{cycle}} = \{v_1, v_2, \dots, v_k\},
\]
where $\{v_k, v_{k+1}\} \notin E(T)$, and the vertex $v_{k+1}$ is the uniquely determined cut vertex. 
\end{lemma}

\noindent
{\bf Proof.}
By Lemma \ref{lem01}, the subgraph $T_{\mathrm{cycle}}$ is a path consisting of a consecutive segment of the cycle. 
The rule that selects the vertex immediately after the endpoint of this segment as the cut vertex is uniquely determined up to a rotation of the cycle labels. 
\qed

\begin{lemma} \label{lem03}
We define a mapping $\varphi : \tau_{n+1} \rightarrow \mathscr{T}_{n}$, 
which specifies the correspondence of vertices and edges as follows. 

\begin{itemize}
\item[$\bullet$]
\textbf{Vertex correspondence}
\begin{itemize}
\item[$\bullet$]
The central vertex $v_c \in W_{n+1}$ corresponds to the hub vertex $v_h \in F_n$. 

\item[$\bullet$]
Each cycle vertex $v_i \in W_{n+1}$ corresponds to the path vertex $u_i \in F_n$. 
\end{itemize}

Here $1 \le i \le n-1$, since the cut vertex $v_n$ is excluded from this correspondence. 
Moreover, the correspondence between cycle edges and path edges is interpreted after reindexing the vertices of the path as in Lemma \ref{lem02}. 

\item[$\bullet$]
\textbf{Edge correspondence}
\begin{itemize}

\item[(i)]
\textit{Cycle edges:}\\
If $T_{\mathrm{cycle}}$ contains the edge 
\(\{v_i, v_{i+1}\}\) for \(1 \le i \le k-1\),
then it corresponds to the path edge \(\{u_i, u_{i+1}\}\).

\item[(ii)]
\textit{General spokes:}\\
If $T_{\mathrm{center}}$ contains a spoke \(\{v_c, v_j\}\) with \(j \ne n\) and \(j\) not being the endpoint of $T_{\mathrm{cycle}}$, then it corresponds to the edge \(\{v_h, u_j\}\). This rule corresponds precisely to spokes at non–cut vertices.

\item[(iii)]
\textit{Spoke at the cut vertex:}\\
If $T_{\mathrm{center}}$ contains the spoke \(\{v_c, v_n\}\), where \(v_n\) is the cut vertex, then it corresponds to the edge \(\{v_h, u_k\}\), where \(k\) is the length of the path in Lemma \ref{lem02}.
\end{itemize}

\end{itemize}

\noindent
\textbf{Remark.}
In the special case where \(T_{\mathrm{center}}=\{v_c\}\), the correspondence given only by rule (i) does not produce a spanning tree.  
Therefore, we add the following rule:

\begin{itemize}
\item[(iv)]
Add the edge \(\{v_h, u_1\}\) to the image. 
\end{itemize}

In what follows, we refer to this special treatment as rule (iv). 
\end{lemma}

\noindent
{\bf Proof.}
For any $T \in \tau_{n+1}$, the number of edges is
\[
|E(T)| = (n+1)-2 = n-1. 
\]
Since the mapping $\varphi$ sends each edge of $T$ to a distinct edge of $\varphi(T)$, we also have $|E(\varphi(T))| = n-1$.
Moreover, every edge of $W_{n+1}$ is mapped either to a hub edge or to a path edge of $F_n$; 
in other words, the image is contained in $F_n$. 

The cycle component $T_{\mathrm{cycle}}$ is mapped to the path
\[
u_1 - u_2 - \cdots - u_k,
\]
and each spoke $\{v_c, v_i\}$ is mapped to an edge incident with $v_h$. 
Therefore the resulting graph is connected and contains exactly $n-1$ edges. 
Hence $\varphi(T)$ is a spanning tree of $F_n$. 
\qed

According to the algorithm described in Lemma \ref{lem03}, we can explicitly verify several examples of the mapping. 
In the following figure, the spanning forest
\[
T = \{v_2v_3,\; v_cv_1,\; v_cv_4\} \in \tau_4
\]
is mapped to the spanning tree
\[
\varphi(T) = \{u_1u_2,\; v_hu_2,\; v_hu_3\} \in \mathscr{T}_3.
\]
Here the cut vertex is $v_4$, and by the edge correspondence (iii), the edge $\{v_c, v_4\}$ is mapped to $\{v_h, u_2\}$. 
This is noted in the figure.
\begin{figure}[H]
\begin{center}
\includegraphics[width=10cm]{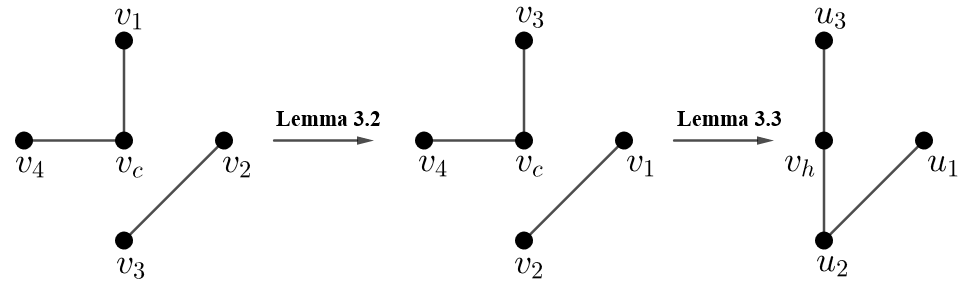}
\end{center}
\end{figure}

In the next example, the spanning forest
\[
T=\{v_1v_2,\; v_2v_3,\; v_3v_4,\; v_c\} \in \tau_4
\]
is mapped to the spanning tree
\[
\varphi(T)=\{u_1u_2,\; u_2u_3,\; v_hu_1\}.
\]
By edge correspondence (iv), if the forest contains an isolated central vertex $v_c$, the additional edge $\{v_h,u_1\}$ is added so that the image becomes a spanning tree. 

\begin{figure}[H]
\begin{center}
\includegraphics[width=10cm]{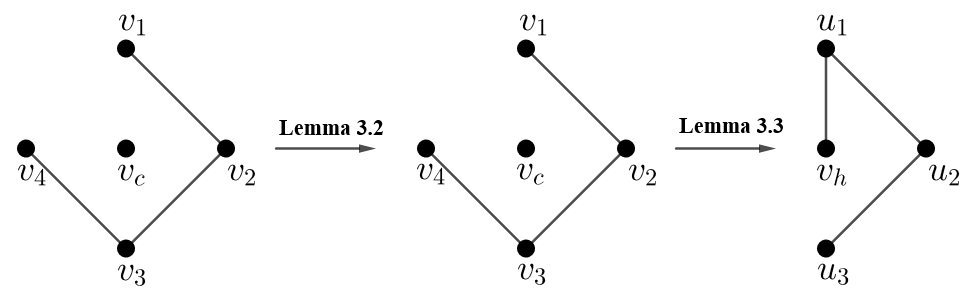}
\end{center}
\end{figure}

\begin{lemma} \label{lem04}
The mapping $\varphi$ is injective. 
Namely, if $\varphi(T)=\varphi(T')$, then $T=T'$. 
\end{lemma}

\noindent
{\bf Proof. }
We show that every $T \in \tau_{n+1}$ can be uniquely reconstructed from its image $S=\varphi(T)$. 
The reconstruction proceeds in the following steps (1) to (4). 

\begin{itemize}
\item[(1)]
{\bf Reconstruction of $T_{\text{cycle}}$. }
For $1 \leq i \leq n-2$, we collect all path edges $\{u_i, u_{i+1}\} \in S$. 
Since $S$ is a spanning tree, these edges form a single continuous path. 
Let this path be $u_1-u_2-\cdots-u_k$. 
This uniquely corresponds to the cycle edges $v_1-v_2-\cdots-v_k$ of $T_{\text{cycle}}$. 

\item[(2)]
{\bf Reconstruction of the cut vertex and its spoke. }
If $\{v_h, u_k\} \in E(S)$, then by rule (iii) of the edge correspondence, 
this edge arises from the spoke $\{v_c, v_n\} \in E(T)$. 
If $\{v_h,u_k\} \notin E(S)$, then $\{v_c,v_n\} \notin E(T)$. 

\item[(3)]
{\bf Reconstruction of the remaining spokes. }
Each edge $\{v_h, u_i\} \in E(S)$ comes from rule (iv), and corresponds to the spoke $\{v_c, v_i\}$ for $1 \leq i \leq n-1$. 

\item[(4)]
{\bf Central vertex is isolated ($T_{\text{center}}=\{v_c\})$. }
If the number of cycle vertices is $k=n$, and the only hub edge in $S$ is $\{v_h, u_1\}$, then this corresponds to rule (iv). 
The reconstructed graph has no spokes, and hence the central vertex $v_c$ is correctly recovered as an isolated vertex. 

\end{itemize}

Thus, by the following steps (1) to (4), the preimage T is uniquely determined from its image $S$. 
Therefore, if $\varphi(T)=\varphi(T')$, then the reconstructed graphs coincide, and hence $T=T'$. 
\qed

\begin{lemma} \label{lem05}
The mapping $\varphi$ is surjective. 
That is, for every $S \in \mathscr{T}_n$, there exists a $T \in \tau_{n+1}$ such that $\varphi(T)=S$. 
\end{lemma}

\noindent
{\bf Proof.}
The argument is similar to the proof of Lemma \ref{lem04}. 
We reconstruct $T$ from $S$ in the following steps (1) to (4). 

\begin{itemize}

\item[(1)]
{\bf Reconstruction of $T_{\text{cycle}}$. }
The path edges $\{u_i, u_{i+1}\} \in S$ form $u_1-u_2-\cdots-u_k$ .
The integer $k$ is the number of vertices in the cycle component. 
Correspondingly, we set $T_{\text{cycle}} = v_1-v_2-\cdots-v_k$ .
By Lemma \ref{lem02}, the cut vertex is $v_{k+1} = v_n$.

\item[(2)]
{\bf Reconstruction of the cut vertex and its spoke.}
I$f \{v_h,u_k\} \in E(S)$, then we include the spoke $\{v_c, v_n\} \in E(T)$. 
If $\{v_h,u_k\} \notin E(S)$, we do not add this edge. 

\item[(3)]
{\bf Reconstruction of the remaining spokes.}
If $\{v_h, u_i\} \in E(S)$, then we include the corresponding spoke $\{v_c, v_i\} \in E(T)$. 

\item[(4)]
{\bf Central vertex isolated.}
If $S$ consists only of the path $u_1-u_2-\cdots-u_{n-1}$ together with the single hub edge $\{v_h,u_1\}$, 
then by rule (iv), no spokes appear in $T$.
Thus the central vertex is correctly reconstructed as isolated.

\end{itemize}

Hence, by the following steps (1) to (4), $T_{\text{cycle}}$ forms a path (and not a cycle),
and the component $T_{\text{center}}$ consisting of $v_c$ and its spokes
forms a star graph. 
Moreover, since no bridges connect the two components, $T$ is a two-component spanning forest. 
This reconstruction procedure is precisely the inverse of rules (i) to (iv)
defining $\varphi$. 
Therefore, for every $S \in \mathscr{T}_n$, there exists a $T \in \tau_{n+1}$ with $\varphi(T)=S$. 
Thus, the mapping $\varphi$ is surjective.
\qed

By the algorithm of Lemma \ref{lem05}, we verify several concrete examples showing that a spanning forest can indeed be recovered from its corresponding spanning tree. 
The following figure illustrates the recovery of the spanning forest
$\tau_4=\{v_1v_2,\, v_cv_3,\, v_cv_4\}$
from the spanning tree
$\mathscr{T}_3=\{u_1u_2,\, v_hu_2,\, v_hu_3\}$ .

\begin{figure}[H]
\begin{center}
\includegraphics[width=12cm]{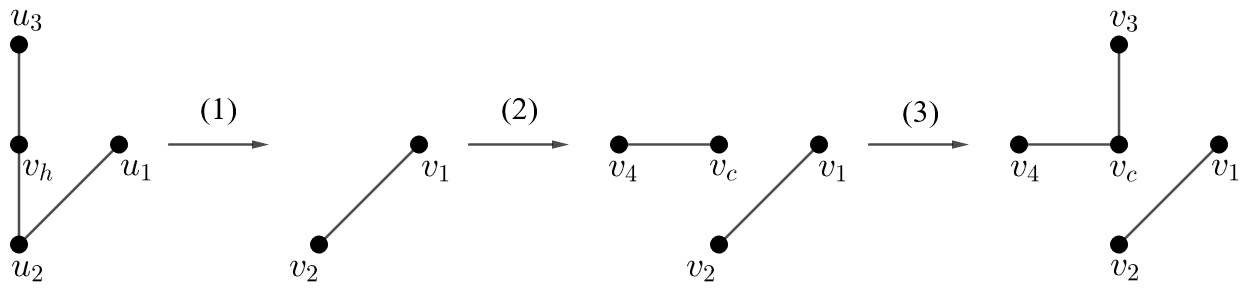}
\end{center}
\end{figure}

The next figure shows the reconstruction of
$\tau_4=\{v_1v_2, v_2v_3, v_3v_4, v_c\}$ from the spanning tree $\mathscr{T}_3=\{u_1u_2, u_2u_3, v_hu_1\}$. 
This example illustrates that when the central vertex is isolated, the reconstruction follows step (4). 

\begin{figure}[H]
\begin{center}
\includegraphics[width=6cm]{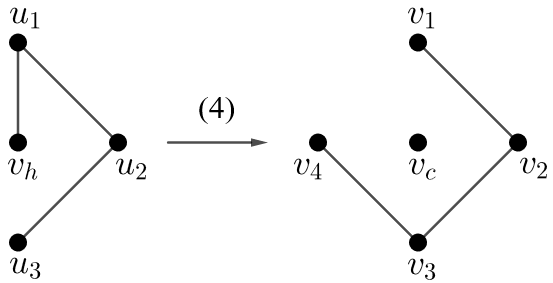}
\end{center}
\end{figure}

From Lemma \ref{lem01} through Lemma \ref{lem05}, the proof of Theorem \ref{thm01} is complete. 

\begin{theorem}[Bogdanowicz, 2008 \cite{P1}] \label{thm03}
Let $F_{n+1}$ be the fan graph consisting of the hub vertex $v_h$ and the path vertices $u_1, u_2, \dots, u_n$. 
The number of spanning trees of $F_{n+1}$ is given by 
\[
\mathscr{T}(F_{n+1}) = f_{2n},
\] 
where $f_{2n}$ denotes the Fibonacci number of index $2n$. 
\end{theorem}

Combining Theorem \ref{thm01} with Theorem \ref{thm03}, we obtain the following corollary. 

\begin{corollary} \label{cor01}
In the wheel graph $W_{n+1}$, let $\tau_{n+1}$ denote the set of two-component spanning forests in which one component contains the central vertex $v_c$ and the other component consists solely of cycle vertices. 
Then the number of such spanning forests satisfies 
\begin{align*}
|\tau_{n+1}|=\mathscr{T}(F_n)=f_{2n-2}. 
\end{align*}
\end{corollary}

In Theorem \ref{thm01}, we established a bijection between 
 • the set $\tau_{n+1}$ of two-component spanning forests of the wheel graph $W_{n+1}$,
where one component contains the central vertex $v_c$ and the other component consists only of cycle vertices, and
 • the set $\mathscr{T}_n$ of spanning trees of the fan graph $F_n$. 

Consequently, by Corollary \ref{cor01}, the number of conditioned spanning forests in the $W_{n+1}$, namely $|\tau_{n+1}|$, coincides with the number of spanning trees of the $F_n$. 

In what follows, we determine the number of two-component spanning forests of $W_{n+1}$ in which two adjacent cycle vertices lie in different components. 
By rotational symmetry of the wheel, we may assume without loss of generality that the two adjacent vertices are $v_1$ and $v_2$, and we compute $F_{W_{n+1}}(v_1 \mid v_2)$. 

\begin{proposition} \label{prop01}
In the $W_{n+1}$, fix the adjacent cycle vertices $v_1$ and $v_2$. 
The number of two-component spanning forests in which $v_1$ and $v_2$ lie in different components is given by
\[
F_{W_{n+1}}(v_1 \mid v_2)=2(f_{2n-1}-1),
\]
where $f_i$ denotes the $i$-th Fibonacci number with $f_0=0$ and $f_1=1$. 
\end{proposition}

\noindent
{\bf Proof. }
Since $v_1$ and $v_2$ must lie in different components, the cycle edge $\{v_1, v_2\}$ cannot be used; otherwise the two vertices would belong to the same component. 
Removing this edge turns the outer cycle into the path $v_x = v_1, v_2, \dots, v_n = v_y$. 
Any spanning forest with two components must place the central vertex $v_c$ in exactly one of the two components. 
Thus the possibilities are: 

\begin{itemize}
\item[(A)]
$v_c$ belongs to the component containing $v_1$, 

\item[(B)]
$v_c$ belongs to the component containing $v_2$. 
\end{itemize}

By symmetry it suffices to consider case (A). 
In this case, the component containing $v_c$ consists of $v_1, v_2, \dots, v_k$ \quad ($1 \leq k \leq n-1$), 
while the remaining vertices $\{v_{k+1}, \dots, v_n\}$ form the other component. 
Since the latter component is a path, its spanning tree is unique.

The component containing $v_c$ has the structure of the fan graph $F_k$, whose number of spanning trees is $\mathscr{T}(F_k)=f_{2k}$. 

Therefore, for fixed $k$ and using the identity
\begin{align*}
F_{W_{n+1}}^{(A)}(v_1 \mid v_2)
&=\sum_{k=1}^{n-1}f_{2k} \\
&=f_2+f_4+f_6+\cdots+f_{2n-2} \\
&=(f_3-f_1)+(f_5-f_3)+(f_7-f_5)+\cdots+(f_{2n-1}-f_{2n-3}) \\
&=f_{2n-1}-1
\end{align*}

we obtain
$F_{W_{n+1}}^{(A)}(v_1 \mid v_2)=f_{2n-1}-1$. 

By symmetry, the same count holds for case (B).
Hence,
\begin{align*}
F_{W_{n+1}}(v_1 \mid v_2)=F_{W_{n+1}}^{(A)}(v_1 \mid v_2)+F_{W_{n+1}}^{(B)}(v_1 \mid v_2)=2(f_{2n-1}-1). 
\end{align*}
\qed

\section{Proof of Theorem \ref{thm02} \label{sec04}}

In this section, we obtain the number of spanning forests of $W_{n+1}$ by using
the effective resistance of the graph. 
Barrett et al. \cite{BarrettEtAl} derived a relation between
two-component spanning forests and effective resistance. 
They proved that the relation among the effective resistance $r_G(u,v)$, 
the number of two-component spanning forests $F_G(u \mid v)$, and the number of spanning trees $\mathscr{T}(G)$ holds, via the all-minors matrix tree theorem of $G$. 

\begin{lemma}[Barrett et al. \cite{BarrettEtAl}] \label{eq:two-forest-resistance}
Let $G$ be a connected graph, and let $u \neq v$ be two vertices of $G$. 
Then the effective resistance $r_G(u,v)$, the number of two-component spanning forests $F_G(u \mid v)$, and the number of spanning trees $\mathscr{T}(G)$ satisfy
\begin{align*}
r_G(u, v)=\frac{F_G(u \mid v)}{\mathscr{T}(G)}.
\end{align*}
Here, $r_G(u,v)$ denotes the effective resistance between $u$ and $v$ in the electrical network obtained by assigning unit resistance to each edge of $G$. 
\end{lemma}

The number of spanning trees of $W_{n+1}$ is well known. 
Myers \cite{myers} and Bogdanowicz \cite{P1} studied the relation between
Fibonacci numbers and Lucas numbers on $W_{n+1}$. 
In particular, the following relation holds for $W_{n+1}$. 

\begin{lemma}[Myers, 1971 \cite{myers}]
\label{lem:wheel-trees}
The number of spanning trees $\mathscr{T}(W_{n+1})$ of $W_{n+1}$ satisfies 
\[
\mathscr{T}(W_{n+1})=\ell_{2n}-2.
\]
Here, $\ell_j$ denotes the $j$th Lucas number, defined by the sequence
$\ell_0=2$, $\ell_1=1$, and $\ell_{j+2}=\ell_{j+1}+\ell_{j}$. 
\end{lemma}

Next, we consider the effective resistance of $W_{n+1}$. 
Bapat and Gupta \cite{bapat-gupta} obtained an explicit formula for the effective resistance $r_G(u, v)$ between any two vertices $u$ and $v$, expressed in terms of Fibonacci numbers, by directly computing the determinant of a principal submatrix of the Laplacian. 

\begin{lemma}[R. Bapat and S. Gupta, 2010 \cite{bapat-gupta}]
\label{eq:bapat}
Let the distance between vertices $u$ and $v$ on the rim cycle be $k=\min\bigl(|u-v|,\; n-|u-v|\bigr)$. 
Then the effective resistance $r_{W_{n+1}}(u, v)$ in $W_{n+1}$ is given by 
\[
r_{W_{n+1}}(u, v)
=\frac{f_{2n}^2}{\,f_{4n}-2f_{2n}\,}\left(2-\frac{f_{4k}}{f_{2k}}\right)+f_{2k}. 
\]
\end{lemma}

Using Lemma \ref{eq:two-forest-resistance} together with Lemma \ref{eq:bapat}, we can express the number of two-component spanning forests of $W_{n+1}$ separating any two vertices $u$ and $v$. 

\noindent
{\bf Proof. }
From Lemma \ref{eq:two-forest-resistance} and Lemma \ref{eq:bapat}, we have
\begin{align*}
F_{W_{n+1}}(u \mid v)
&=r_{W_{n+1}}(u, v) \cdot \mathscr{T}(W_{n+1}) \\
&=\left\{\frac{f_{2n}^2}{f_{4n}-2f_{2n}}\left(2-\frac{f_{4k}}{f_{2k}}\right)+f_{2k} \right\} \cdot \mathscr{T}(W_{n+1})
\end{align*}

Next, applying the Fibonacci and Lucas identity $f_{4n}=f_{2n}\ell_{2n}$ together with Lemma \ref{lem:wheel-trees}, we obtain
\[
f_{4n}-2f_{2n}
=f_{2n}(\ell_{2n}-2)
=f_{2n}\,\mathscr{T}(W_{n+1}).
\]

Proceeding with the calculation further, we derive
\begin{align*}
F_{W_{n+1}}(u \mid v)
&=\left\{\frac{f_{2n}^2}{f_{2n} \cdot \mathscr{T}(W_{n+1})}\left(2-\frac{f_{4k}}{f_{2k}}\right)+f_{2k} \right\} \cdot \mathscr{T}(W_{n+1}) \\
&=f_{2n}\left(2-\frac{f_{4k}}{f_{2k}}\right)+f_{2k} \mathscr{T}(W_{n+1}) \\
&=f_{2n}\left(\frac{2f_{2k}-f_{2k}\ell_{2k}}{f_{2k}}\right)+f_{2k} \mathscr{T}(W_{n+1}) \\
&=f_{2n}(2-\ell_{2k})+f_{2k} \mathscr{T}(W_{n+1}) \\
&=f_{2k} \mathscr{T}(W_{n+1})-f_{2n} \mathscr{T}(W_{k+1}) \\
&=f_{2k}(\ell_{2n}-2)-f_{2n}(\ell_{2k}-2)
\end{align*}

Therefore, the proof of Theorem \ref{thm02} is complete. 
\qed

Using Theorem \ref{thm02}, we obtain the numbers of the three types of conditioned two-component spanning forests: 
$F_{W_{n+1}}(v_1\mid v_2)$, $F_{W_{n+1}}(v_1\mid v_3)$, $F_{W_{n+1}}(v_1\mid v_c)$.

\begin{corollary} \label{corollary41}
The number of two-component spanning forests of $W_{n+1}$ in which two adjacent vertices, such as $v_1$ and $v_2$, lie in different components is given by
\[
F_{W_{n+1}}(v_1\mid v_2)=2(f_{2n-1}-1),
\]
where $f_{2n-1}$ denotes the Fibonacci number of index $2n-1$. 
\end{corollary}

\noindent
{\bf Proof.}
In this case, we substitute the cycle distance $k=1$, together with the Fibonacci number $f_2=1$ and the Lucas number $\ell_2=3$, into
Theorem \ref{thm02}: 
\begin{align*}
F_{W_{n+1}}(v_1\mid v_2)
&=f_2(\ell_{2n}-2)-f_{2n}(\ell_2-2) \\
&=\ell_{2n}-2-f_{2n} \\
&=2f_{2n-1}-2 \\
&=2(f_{2n-1}-1). 
\end{align*}
The third equality uses the identity $\ell_{2n}-f_{2n}=2f_{2n-1}$. 
\qed

This result agrees with Proposition \ref{prop01}. 

\begin{corollary} \label{corollary42}
The number of two-component spanning forests of $W_{n+1}$ in which two vertices at cycle distance $k=2$, such as $v_1$ and $v_3$, lie in different components is
\[
F_{W_{n+1}}(v_1\mid v_3)=2(\ell_{2n-2}-3)=2(\mathscr{T}(W_{n})-1), 
\]
where $\ell_{2n-2}$ is the Lucas number of index $2n-2$, and $\mathscr{T}(W_{n})$ denotes the number of spanning trees of $W_{n}$. 
\end{corollary}

\noindent
{\bf Proof. }
In this case, we substitute the cycle distance $k=2$, together with the Fibonacci number $f_4=3$ and the Lucas number $\ell_4=7$ into Theorem \ref{thm02}, we obtain
\begin{align*}
F_{W_{n+1}}(v_1\mid v_3)
&=f_4(\ell_{2n}-2)-f_{2n}(\ell_4-2) \\
&=3(\ell_{2n}-2)-f_{2n}(7-2) \\
&=3\ell_{2n}-5f_{2n}-6 \\
&=2\ell_{2n-2}-6 \\
&=2(\ell_{2n-2}-3)
=2(\mathscr{T}(W_{n})-1).
\end{align*}
The fourth equality uses the identity $\ell_{2n}-5f_{2n} = 2\ell_{2n-2}$. 
\qed

\begin{corollary} \label{corollary43}
The number of two-component spanning forests of $W_{n+1}$ in which $v_1$ and the central vertex $v_c$ lie in different components is
\[
F_{W_{n+1}}(v_1\mid v_c)=f_{2n},
\]
where $f_{2n}$ denotes the Fibonacci number of index $2n$. 
\end{corollary}

\noindent
{\bf Proof.}
The effective resistance between the cycle vertex $v_1$ and the central vertex $v_c$ is given by Bapat and Gupta \cite{bapat-gupta}:
\[
r(v_1, v_c)=\frac{f_{2n}^{2}}{f_{4n}-2f_{2n}}. 
\]
Using $f_{4n}-2f_{2n}=f_{2n} \mathscr{T}(W_{n+1})$, we obtain
\[
r(v_1, v_c)=\frac{f_{2n}}{\mathscr{T}(W_{n+1})}. 
\]
Hence, by Theorem \ref{eq:two-forest-resistance},
\[
F_{W_{n+1}}(v_1\mid v_c)
=\mathscr{T}(W_{n+1}) \cdot r(v_1, v_c)
=f_{2n}.
\]
\qed

From Corollary \ref{corollary41} through Corollary\ref{corollary43}, we complete the proof of Proposition \ref{prop1a}. 

\section{Conclusion} \label{sec05}

In this paper, we studied the structure of conditioned two-component spanning forests on $W_{n+1}$ using both combinatorial and analytic methods.

First, we constructed a bijection between the set of two-component spanning forests $\tau_{n+1}$ and the set of spanning trees
$\mathscr{T}_n$ of the fan graph $F_n$. 
This correspondence is obtained by the edge operations, and we proved that the bijection holds. 
As a result, we obtained $|\tau_{n+1}|=\mathscr{T}(F_n)=f_{2n-2}$ and found that the Fibonacci structure is preserved between the distinct structures of $W_{n+1}$ and $F_n$. 

Second, by applying the all minors matrix tree theorem of Barrett et al. \cite{BarrettEtAl} and the explicit effective resistance formulas
of Bapat and Gupta \cite{bapat-gupta}, we obtained the number of two-component spanning forests in which any two vertices $u$ and $v$ belong to different components.
Using this result, we derived the three explicit formulas for $F_{W_{n+1}}(v_1\mid v_2)$, $F_{W_{n+1}}(v_1\mid v_3)$, and $F_{W_{n+1}}(v_1\mid v_c)$, all expressed in terms of Fibonacci and Lucas numbers. 

Our results were not previously listed in the OEIS; we submitted these sequences, and they have been published as A391263. 
Moreover, they have also been published as new combinatorial interpretations in A001906 and A027941. 
This work integrates two complementary viewpoints: a constructive understanding via a bijection and an analytic understanding via
effective resistance. It provides an example demonstrating how the spanning structure of a graph is related to the structure of integer sequences. 
As a direction for future work, it would be interesting to provide combinatorial interpretations of two-component spanning forests of other
graph families and to investigate which graphs admit such bijections.


\end{document}